\documentclass[12pt,reqno]{amsart}
\usepackage{amssymb}
\usepackage{txfonts}
\usepackage{amsmath}

 \newtheorem{thm}{Theorem}[section]

 \theoremstyle{definition}
 \newtheorem{defn}[thm]{Definition}
 \theoremstyle{remark}
 \newtheorem{rem}[thm]{Remark}
 \numberwithin{equation}{section}

\begin{document}

\title[]
{An asymptotic property of Huisken's functional on  minimal submanifolds of Euclidean space}

\author{Liang Cheng$^1$}

\dedicatory{}
\date{}


\thanks{$^1$Research partially supported by the Natural Science
Foundation of China 11201164 and 11171126}

\address{Liang Cheng, School of Mathematics and Statistics, Huazhong Normal University,
Wuhan, 430079, P.R. CHINA}

\email{math.chengliang@gmail.com}

 \maketitle

\begin{abstract}
In this short note, we study the asymptotic property of Huisken's functional for mean curvature flow on the minimal submanifolds of Euclidean space. We prove that the limit of Huisken's functional equals to the extrinsic asymptotic volume ratio on the minimal submanifold of Euclidean space.
\end{abstract}

\section{Introduction}
Huisken introduced his entropy in \cite{H1} which becomes one of powerful tools for studying mean curvature flow.
Huisken's entropy enjoys very nice analytic and geometric properties,
including in particular the monotonicity of the entropy. These
properties can be used, as demonstrated by Huisken \cite{H1}, Ilmanen \cite{I1}, White \cite{W1} to show the singularities of the mean curvature flow can be modelled by self-shrinking solutions of the flow.

Let $X:M^n\to\mathbb{R}^{n+m}$ be the complete minimal immersed submanifold with $M=X(M^n)$.
Recall the Huisken's entropy is defined as the integral of backward heat kernel:
\begin{equation}\label{Huisken}
\int_{M}[4\pi (t_0-t)] ^{-\frac{n}{2}}e^{-\frac{|x-y_0|^2}{4(t_0-t)}}d\mu(x),\ t\in (-\infty,t_0),\ y_0\in\mathbb{R}^{n+m}.
\end{equation}
Huisken \cite{H1} proved his entropy (\ref{Huisken}) is monotone non-increasing in $t$ under the mean curvature flow
$
\frac{\partial X}{\partial t}=-H\nu,
$
where $H$ is the mean curvature of $X$ and $\nu$ is the outer unit
normal to $M$. Moreover, (\ref{Huisken}) is invariant in $t$ if and only if $M$ is the cone in $\mathbb{R}^{n+m}$ (see \cite{H1}) .
Let $\tau=t_0-t$. Denote
\begin{equation}\label{Huisken_1}
H_{y_0,\tau}(M)=\int_{M}(4\pi \tau) ^{-\frac{n}{2}}e^{-\frac{|x-y_0|^2}{4\tau}}d\mu(x),\ \tau\in (0,+\infty),\ y_0\in\mathbb{R}^{n+m}.
\end{equation}
Then $H_{y_0,\tau}(M)$ is non-decreasing in $\tau$ under the backward mean curvature flow
$
\frac{\partial X}{\partial \tau}=H\nu
$.

Though the mean curvature flow is static on minimal submanifolds,  we remark that Huisken's functional $H_{y_0,\tau}(M)$ is still strictly increasing in $\tau$ on the minimal submanifold $M$ unless $M$ is the minimal cone.
So it is interesting to understand what's the geometric meaning of the limit of Huisken's functional $H_{y_0,\tau}(M)$ on minimal submanifolds of Euclidean space.

 Note that $\frac{Vol(B(y_o,r)\cap M)}{r^n}$ is an non-decreasing function of $r$ by the monotonicity formula for the minimal submanifolds of Euclidean space (see Proposition 1.8 in \cite{CM2}), where
$B(y_0, r)$ is the ball in Euclidean space centered at $y_0$ with radius $r$. Hence we use the following
\begin{defn}
Let $X:M^n\to\mathbb{R}^{n+m}$ be the complete minimal immersed submanifold with $M=X(M^n)$.
The extrinsic asymptotic volume ratio ($EAVR$) is defined as the following
\begin{equation}
EAVR(y_0)=\lim\limits_{r\to\infty}\frac{Vol(B(y_0,r)\cap M)}{\omega_n r^n},
\end{equation}
where
$\omega_n$ is the volume of unit $n$-ball in
$\mathbb{R}^{n+m}$.
\end{defn}
There are many examples of minimal submanifolds that $EAVR(y_0)$ is finite. It is clear that an nontrivial example which has finite $EAVR(y_0)$ is the catenoid in $\mathbb{R}^3$. Moreover, the minimal graphs have the finite $EAVR(y_0)$ (see \cite{CM2}). According to the curvature estimate of Osserman \cite{OS}, any complete minimal
surfaces in $\mathbb{R}^3$ with finite total curvature has finite $EAVR(y_0)$. We also make the following observation. Consider tangent cones at infinity along the following way and define, for $r_j>0$,
$$
M_{r_j}=\frac{1}{r_j}\{B(y_0,r_j)\cap M-\{y_0\}\}
$$
and
$$
\Sigma_{r_j}=\frac{1}{r_j}\{\partial B(y_0,r_j)\cap M-\{y_0\}\}.
$$
Then $M_{r_j}\subset B(0,1)\subset \mathbb{R}^{n+m}$, $\Sigma_{r_j}\subset \partial B(0,1)\subset S^{n+m-1}$ and
$$EAVR(y_0)=\lim\limits_{r_j\to\infty}Vol(M_{r_j})=\lim\limits_{r_j\to\infty}Area(\Sigma_{r_j}).$$
For the minimal surface
with finite $EAVR(y_0)$, one can define the asymptotic minimal
varifold limit at infinity by taking any sequence $r_j\to\infty $.

We show that the limit of Huisken's functional equals to  $EAVR(y_0)$ on minimal submanifolds of Euclidean space.
\begin{thm}\label{main1}
Let $X:M^n\to\mathbb{R}^{n+m}$ be the complete immersed minimal submanifold with $M=X(M^n)$. Assume that $H_{y_0, \tau}(M)$ is well defined for all $\tau>0$.  Then $\lim\limits_{ \tau \to\infty}H_{y_0, \tau}(M)=EAVR(y_0)$,
\end{thm}
\begin{rem}
Clearly, $H_{y_0,\tau} (M)$ is well defined for all $\tau>0$ on a complete
submanifold in $\mathbb{R}^{n+m}$ with bounded second fundamental form.
\end{rem}
\begin{rem}
We should mention that there is a similar result to Theorem \ref{main1} that $\lim\limits_{\tau\to\infty}\int_{M} (4\pi\tau)^{-\frac{n}{2}}e^{-\frac{d_g(p,y)^2}{4\tau}}dvol_{g}(y)$
equals to the intrinsic asymptotic volume ratio (i.e. $\lim\limits_{r\to\infty}\frac{\text{Vol}B_g(p,r)}{\omega_n r^n}$) on the $n$-dimensional  Riemannian manifold $(M^n,g)$ with nonnegative Ricci curvature (see \cite{CCGGI}), where $d_g$ is the distance on $(M^n,g)$ and $B_g(p,r)$ is the ball on $(M^n,g)$ centered at $p$ with radius $r$.
\end{rem}

\section{Proof of Theorem \ref{main1}}
In this section, we give the proof of Theorem \ref{main1}.

\textbf{Proof of Theorem \ref{main1}.}
Let $A(s)= Area(\partial B(y_0,s)\cap M)$. By the monotonicity formula of the minimal submanifolds of Euclidean space (see Proposition
1.8 in \cite{CM2}), we have
$$
\frac{Vol(B(y_0,s)\cap M)}{s^n}\geq \frac{Vol(B(y_0,r)\cap M)}{r^n}
$$
for any $s\geq r$.
Let $f(s)=Vol(B(y_0,s)\cap M)$ and $g(s)=s^n$. Then
$$
\frac{f(s)}{g(s)}\geq \frac{f(r)}{g(r)}
$$
for any $s\geq r$.  It follows that
$$
\frac{\frac{f(s)-f(r_1)}{s-r_1}}{\frac{g(s)-g(r_1)}{s-r_1}}\geq \frac{f(r)}{g(r)}
$$
for any $s\geq r_1\geq r$. Let $r_1\to s$, we get
$$
\frac{f'(s)}{g'(s)}\geq \frac{f(r)}{g(r)}.
$$
for any $s\geq r$. Since $f(s)=\int^s_0A(s)ds$, we have
$$
\frac{A(s)}{ns^{n-1}}\geq \frac{Vol(B(y_0,r)\cap M)}{r^n}
$$
for any $s\geq r$. Similarly, we can prove that
\begin{equation}\label{eq_33}
\frac{A(s)}{ns^{n-1}}\leq \lim\limits_{r\to\infty} \frac{Vol(B(y_0,r)\cap M)}{r^n},
\end{equation}
 for any $s>0$.
It follows that
\begin{align*}
H_{y_0,\tau}(M)&=\int_{ B(y_0,r)}(4\pi \tau) ^{-\frac{n}{2}}e^{-\frac{|x-y_0|^2}{4\tau}}d\mu(x)+\int_{ M\backslash B(y_0,r)}(4\pi \tau) ^{-\frac{n}{2}}e^{-\frac{|x-y_0|^2}{4\tau}}d\mu(x)\\
&= \int_{ B(y_0,r)}(4\pi \tau) ^{-\frac{n}{2}}e^{-\frac{|x-y_0|^2}{4\tau}}d\mu(x)+\int^{\infty}_r (4\pi \tau) ^{-\frac{n}{2}}e^{-\frac{s^2}{4\tau}}A(s)ds\\
&\geq (4\pi \tau) ^{-\frac{n}{2}}e^{-\frac{r^2}{4\tau}}Vol(B(y_0,r)\cap M)+\frac{n}{2}\pi^{-\frac{n}{2}}\frac{Vol(B(y_0,r)\cap M)}{r^n}\int^{\infty}_{\frac{r^2}{4\tau}} e^{-\eta}\eta^{\frac{n-2}{2}} d\eta\\
&= \frac{Vol(B(y_0,r)\cap M)}{r^n}((4\pi \tau) ^{-\frac{n}{2}}e^{-\frac{r^2}{4\tau}}r^n+\frac{n}{2}\pi^{-\frac{n}{2}}\int^{\infty}_{\frac{r^2}{4\tau}} e^{-\eta}\eta^{\frac{n-2}{2}} d\eta).
\end{align*}
Since $\omega_n\frac{n}{2}\pi^{-\frac{n}{2}}\int^{\infty}_{0} e^{-\eta}\eta^{\frac{n-2}{2}} d\eta=\int_{\mathbb{R}^n}\pi^{-\frac{n}{2}}e^{-|x|^2}d\mu(x)=1$, we have
$\lim\limits_{\tau\to\infty}H_{y_0,\tau}(M)\geq \frac{Vol(B(y_0,r)\cap M)}{\omega_nr^n}$ for any $r>0$. Hence
\begin{equation}\label{eq_31}
\lim\limits_{\tau\to\infty}H_{y_0,\tau}(M)\geq \lim\limits_{r\to\infty}\frac{Vol(B(y_0,r)\cap M)}{\omega_nr^n}.
 \end{equation}
By (\ref{eq_33}), we have
\begin{align*}
H_{y_0,\tau}(M)&=\int^{\infty}_0 (4\pi \tau) ^{-\frac{n}{2}}e^{-\frac{s^2}{4\tau}}A(s)ds\\
&\leq n\lim\limits_{r\to\infty} \frac{Vol(B(y_0,r)\cap M)}{r^n}\int^{\infty}_0 (4\pi \tau) ^{-\frac{n}{2}}e^{-\frac{s^2}{4\tau}}s^{n-1}ds\\
&\leq \lim\limits_{r\to\infty} \frac{Vol(B(y_0,r)\cap M)}{\omega_nr^n}.
\end{align*}
It follows that
\begin{align}\label{eq_32}
\lim\limits_{\tau\to\infty}H_{y_0,\tau}(M)\leq \lim\limits_{r\to\infty} \frac{Vol(B(y_0,r)\cap M)}{\omega_nr^n}.
\end{align}
Combining with (\ref{eq_31}) and (\ref{eq_32}), we get
\begin{align*}
\lim\limits_{\tau\to\infty}H_{y_0,\tau}(M)= \lim\limits_{r\to\infty} \frac{Vol(B(y_0,r)\cap M)}{\omega_nr^n}.
\end{align*}
$\Box$

\thanks{\textbf{Acknowledgement}: The author would like to thank the unknown referee for his/her careful reading of the earlier version of this paper and many fruitful comments and suggestions, which improve the quality of this paper.}

\end{document}